\documentclass{amsart}
\pdfoutput=1
\usepackage{caption}
\usepackage{subcaption}
\usepackage{amssymb}
\usepackage{graphicx}
\usepackage{wrapfig}

\usepackage{amssymb}
\usepackage{graphicx}

\author{G. Fejes T\'oth}
\address{Alfr\'ed R\'enyi Institute of Mathematics,
Re\'altanoda u. 13-15., H-1053, Budapest, Hungary}
\email{gfejes@renyi.hu}

\title{Finite variations on the isoperimetric problem}
\thanks{The English translation of the book ``Lagerungen in der Ebene,
auf der Kugel und im Raum" by L\'aszl\'o Fejes T\'oth will be
published by Springer in the book series Grundlehren der
mathematischen Wissenschaften under the title
``Lagerungen---Arrangements in the Plane, on the Sphere and
in Space". Besides detailed notes to the original text the
English edition contains eight self-contained new chapters
surveying topics related to the subject of the book but not
contained in it. This is a preprint of one of the new chapters.}

\begin{document}

\begin{abstract}
The isoperimetric problem asks for
the maximum area of a region of given perimeter. It is natural to consider
other measurements of a region, such as the diameter and width,
and ask for the extreme value of one when another
is fixed. The solution of these problems is known if the competing regions
are general convex disks, however several of these problems are still open
if the competing regions are polygons with at most a given number of sides.
The present work surveys these problems.

\end{abstract}

\maketitle

Let $a$, $p$, $w$, and $d$ denote the area, perimeter, width and diameter
of a convex disk. Fixing one of these four quantities, what is the infimum
and the supremum of another one of them? Of course, fixing one quantity and
asking for the supremum of another one is equivalent to the problem of
fixing the second quantity and asking for the infimum of the first one.
The solution of one half of the twelve problems arising this way is obvious:
The answer is either zero or infinity. In the case of the six meaningful
problems, we can ask for minima and maxima and for the convex disks
attaining the optimum.

The isoperimetric problem asks for the convex disk of maximum area with
given perimeter. Its solution is the circle, which was known already in
Ancient Greece, although a mathematically rigorous proof was obtained
only in the 19th century. The solution of the problem on the maximum
width for a given diameter is obvious: The optimal sets are convex disks of
constant width. Convex disks of constant width were also characterized by
{\sc{Blaschke}} \cite{Blaschke} as those among domains of a given width
that have minimum perimeter, and by {\sc Rosenthal} and {\sc Sz{\'a}sz}
\cite{RosenthalSzasz} as those among  with a given diameter, that have
maximum perimeter. {\sc{Bieberbach}} \cite{Bieberbach} proved that among
all domains of a fixed diameter the one of maximum area is the circular
disk, and {\sc{P\'al}} \cite{Pal21} proved that the minimum area of a convex
disk with given width is attained by the regular triangle.

An interesting area of research is to consider these optimum problems
restricted to polygons with at most a given number of sides. The
isoperimetric problem for $n$-gons was solved centuries ago by Zenodorus
(see \cite{Blasjo}): Implicitly assuming the existence of a solution, he
proved that a regular $n$-gon has greater area than all other $n$-gons
with the same perimeter. This is expressed in the inequality
$$p^2\ge4na\tan\frac{\pi}{n}.$$
The theorem of P\'al solves the problem of minimum area for a given width.
Concerning the remaining four problems only partial results are known.

{\sc Reinhardt} \cite{Reinhardt} considered the problem of maximizing the perimeter
of a convex $n$-gon with a given diameter and proved that the perimeter $p$ of an $n$-gon
of diameter $d$ satisfies the inequality
$$p\le2n\sin\frac{\pi}{2n}d.$$
Equality is attained here if and only if $n$ is not a power of 2. This result, together
with the characterization of the case of equality, was rediscovered by {\sc{Larman}} and
{\sc{Tamvakis}} \cite{LarmanTamvakis}, {\sc{Datta}} \cite{Datta} and {\sc{A.~Bezdek}}
and {\sc{Fodor}} \cite{BezdekAFodor00}.

To describe the polygons for which equality is attained, we start with a convex polygon with an
odd number of sides such that each vertex is at distance $d$ from the endpoints of the opposite
side. Replacing each side by a circular arc of radius $d$ centered at the opposite vertex we
obtain a Reuleaux polygon. If $n$ is not a power of 2, the $n$-gons of diameter $d$ with
perimeter $2n\sin(\pi/2n)$ are inscribed in a Reuleaux polygon in such a way that every vertex
of the Reuleaux polygon is a vertex of the polygon, and all sides of the polygon are of equal
length. Such polygons are called by {\sc{Audet}}, {\sc{Hansen}} and {\sc{Messine}}
\cite{AudetHansenMessine09a} {\it clipped Reuleaux polygons}, while {\sc{Mossinghoff}}
\cite{Mossinghoff11} uses the term {\it{Reinhardt polygons}} for them.

\medskip
\centerline {\immediate\pdfximage width4cm
{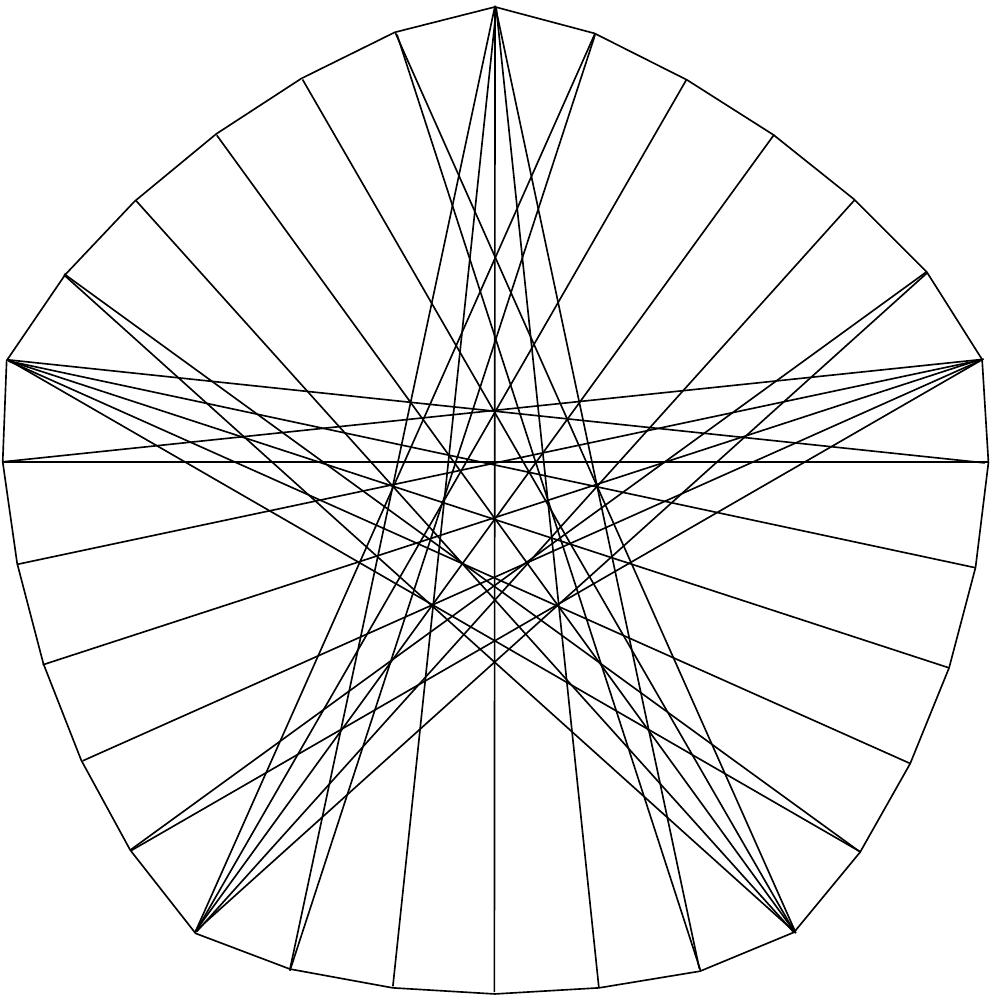}\pdfrefximage \pdflastximage\hskip.5truecm
\immediate\pdfximage width4cm
{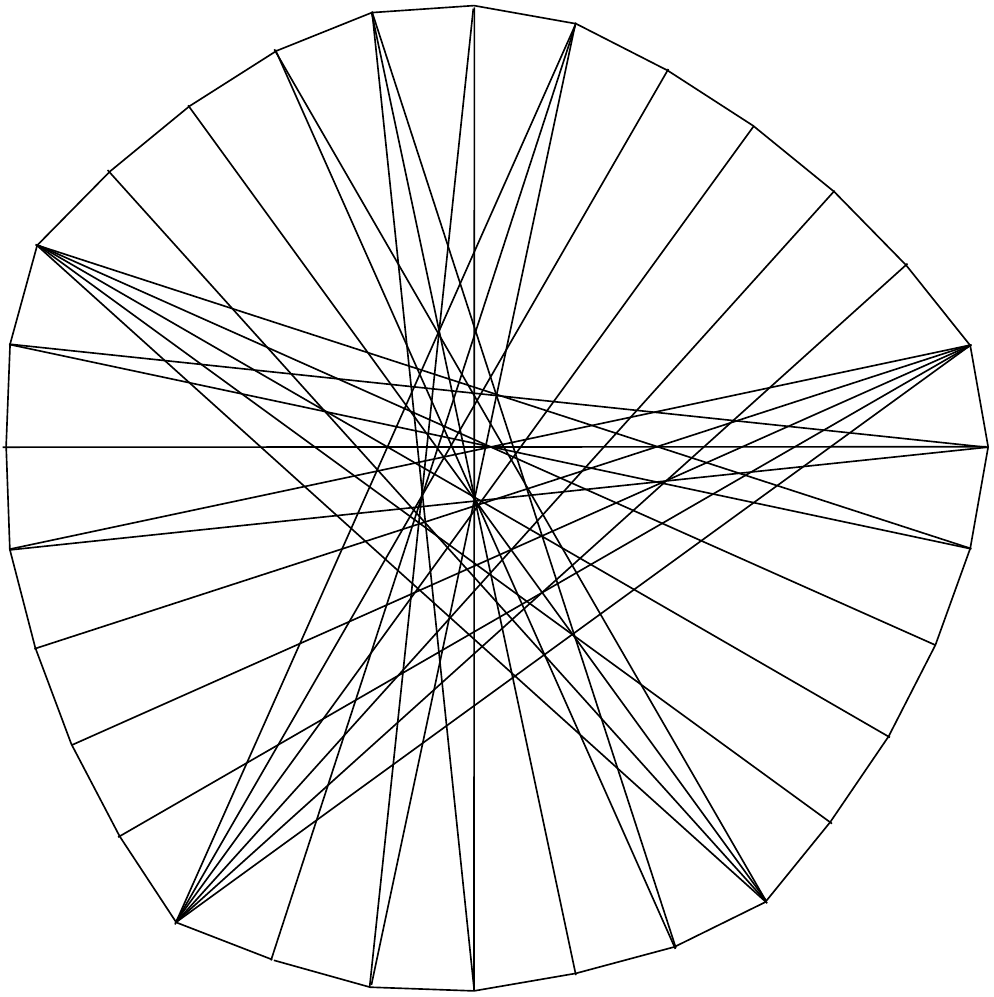}\pdfrefximage \pdflastximage}
\smallskip{\centerline{Figure~1}
\medskip

Consider a Reuleaux polygon with $m$ vertices. Its diagonals form an $m$-gon which for
$m>3$ is a star polygon. The sum of the angles of this polygon is $\pi$, so in order that
it can accommodate a clipped Reuleaux polygon its angles must be integer multiples of $\pi/n$.
The clipped Reuleaux polygons were studied by {\sc{Gashkov}} \cite{Gashkov07,Gashkov13},
{\sc{Mossinghoff}} \cite{Mossinghoff11} and {\sc{Hare}} and  {\sc{Mossinghoff}}
\cite{HareMossinghoff13,HareMossinghoff19}. For a given $n$, there are clipped Reuleaux
$n$-gons with $k$-fold rotational symmetry for some divisor $k$ of $n$. Besides these,
called {\it{periodic}} by Mossinghoff, there may be some others, called {\it{sporadic}}.
The latter name is misleading, since it turned out that the sporadic clipped Reuleaux
polygons outnumber the periodic ones for almost all $n$. Figure~1 shows a regular and
a sporadic clipped Reuleaux polygon with 30 vertices.

Finding the maximum perimeter of an $n$-gon of given diameter when $n$ is a power of 2
is difficult. Only the cases of the quadrangle and octagon are solved. The best quadrangle was
determined by {\sc{Tamvakis}} \cite{Tamvakis} and rediscovered by {\sc{Datta}} \cite{Datta}. The
octagon's case was settled by {\sc{Audet, Hansen}} and {\sc{Messine}} \cite{AudetHansenMessine07a}.
Tamvakis described a sequence of unit-diameter $n$-gons for $n=2^k$ whose perimeter exceeds that of
the regular $n$-gon, and differs from the upper bound $2n\sin\frac{\pi}{2n}$ by $O(n^{-4})$.
By improved constructions the difference from the upper bound was reduced to $O(n^{-5})$ by
{\sc{Mossinghoff}} \cite{Mossinghoff06a} and lately to $O(n^{-6})$ by {\sc{Bingane}} \cite{Bingane21a}.

Combining the inequality $p\le2n\sin\frac{\pi}{2n}d$ with the isoperimertic inequality
$p^2\ge4n\tan\frac{\pi}{n}a$, {\sc Reinhardt} \cite{Reinhardt} obtained the inequality
$$a\le\frac{n}{2}\cos\frac{\pi}{n}\tan\frac{\pi}{2n}d^2$$
with equality only for odd $n$ and regular $n$-gons. Thus, for odd $n$, among all $n$-gons
of a given diameter the regular one has maximum area. Alternative proofs were given by
{\sc{Lenz}} \cite{Lenz56a}, {\sc{Griffiths}} and {\sc{Culpin}} \cite{GriffithsCulpin}, and
{\sc{Gashkov}} \cite{Gashkov85}.

Reinhardt proved that for even $n\ge6$, the optimal $n$-gon is never regular. Alternative
proofs were given by {\sc{Sch\"affer}} \cite{Schaffer}, {\sc{Audet, Hansen}} and {\sc{Messine}}
\cite{AudetHansenMessine08} and {\sc{Mossinghoff}} \cite{Mossinghoff06a}. The latter author
constructed a sequence of $n$-gons with unit diameter for even $n$ whose area exceeds the
area of the unit-diameter regular $n$-gon by $O(n^{-2}$), and whose area differs from the
maximum area of such $n$-gons by a term of at most $O(n^{-3})$. {\sc{Bingane}}
\cite{Bingane21b} improved {\sc{Mossighoff}}'s construction without improving
on the order of difference from the maximum area.

The maximum area of a quadrangle of diameter $d$ is $d^2/2$. The diagonals of the optimal
quadrangles are perpendicular and have length $d$. The case
of the hexagon was solved by {\sc{Graham}} \cite{Graham}. He confirmed the conjecture of
{\sc{Bieri}} \cite{Bieri} that the non-regular hexagon shown in Figure 2 is the unique
optimal solution. He also formulated a conjecture for all even $n\ge6$, stating that for
such $n$, every optimal $n$-gon's diameter graph consists of an ($n-1$)-cycle with one
additional edge emanating from one of the cycle's vertices. Note that the conjecture
leaves the geometric realization of the best polygon undetermined, subject to an
optimization problem. Graham's conjecture was confirmed for $n=8$ by {\sc{Audet, Hansen,
Messine}}, and {\sc{Xiong}} \cite{AudetHansenMessineXiong}, who also solved the
corresponding optimization problem, thus determining the best octagon. Subsequently,
Graham's conjecture was confirmed in general by {\sc Foster} and {\sc{Szabo}}
\cite{FosterSzabo}. The corresponding optimal polygons for $n=10$ and $n=12$ were
determined by {\sc{Henrion}} and {\sc{Messine}} \cite{HenrionMessine}.

\medskip
\centerline {\immediate\pdfximage width4cm
{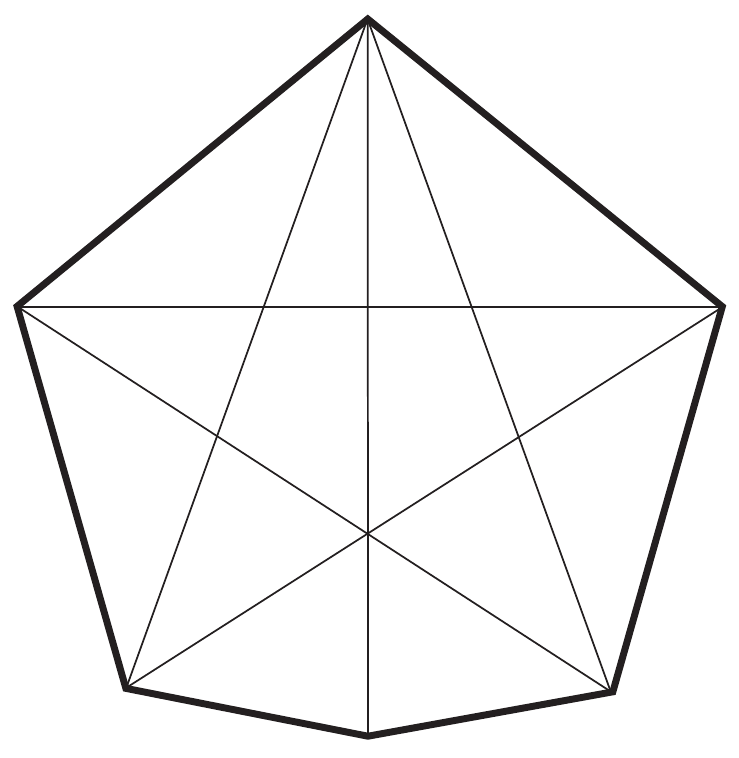}\pdfrefximage \pdflastximage}
\smallskip{\centerline{Figure~2}}
\medskip

{\sc{Graham}} \cite{Graham} also asked for the solution of the higher-dimensional analogue of
the problem: Which convex $d$-polytope with $n$ vertices and unit diameter has the largest volume?
For $n=d+1$ the solution is the regular simplex. {\sc{Kind}} and {\sc{Kleinschmidt}}
\cite{KindKleinschmidt76} solved the problem for $n=d+2$ and described all the extremal
polytopes. The case $n=d+3$ was attacked by {\sc{Klein}} and {\sc{Wessler}} \cite{KleinWessler03},
however their proof turned out to be incomplete (cf. {\sc{Klein}} and {\sc{Wessler}}
\cite{KleinWessler05}), thus, this case is still open.

In the Russian journal for high school students Quant, {\sc{Gashkov}} \cite{Gashkov85} wrote
a small article about the isoperimetrtic problem and its relatives. There he gave a proof of
Reinhard's theorem based on central symmetrization. He used the following facts about
a convex $n$-gon $P$ and its central symmetric image $P^*=\frac{1}{2}(P-P)$: $P^*$ is
a convex polygon with at most $m\le2n$ sides and the same width $w$, diameter $d$ and
perimeter $p$ as $P$. The inradius of $P^*$ is at least $w$ and the circumradius of
$P^*$ is at least $d$. It follows that
$$2m\sin\frac{\pi}{2m}d\ge{p}\ge{m}\tan\frac{\pi}{m}w.$$
Since $m\le2n$, using the monotonicity of the functions $x\sin\frac{1}{x}$ and
$x\tan\frac{1}{x}$, we get, on one hand, Reinhardt's inequality for the maximum
perimeter of an $n$-gon with given diameter, and on the other hand, the new inequality
$$p\ge{2n}\tan\frac{\pi}{2n}w.$$
The combination of these inequalities yields the inequality
$$w\le\cos\frac{\pi}{2n}d$$
between the width and diameter of a convex $n$-gon. Equality is attained in these
inequalities for every $n\ge3$ that has an odd factor by a clipped Reuleaux polygon.

Gashkov's article remained unnoticed. The last inequality was rediscovered by
{\sc{A.~Bezdek}} and {\sc{Fodor}} \cite{BezdekAFodor00}, and the inequality between
perimeter and width by {\sc Audet}, {\sc Hansen} and {\sc Messine}
\cite{AudetHansenMessine09a}. These authors also solved the case of the quadrangle for
both problems. The octagon of a given diameter with maximum width was determined by
{\sc{Audet, Hansen, Messine}} and {\sc{Ninin}} \cite{AudetHansenMessineNinin}.

Motivated by a question of Erd\H{o}s, {\sc{Vincze}} \cite{Vincze} studied the problem of finding
the maximum perimeter of an equilateral $n$-gon with given diameter. He solved the problem if $n$ is
not a power of 2. Of course, this case is an immediate consequence of Reinhardt's theorem. However,
Vincze's argument works without the assumption that the sides have equal length, so it yields an
alternative proof of Reinhardt's theorem. The case that $n$ is a power of $2>4$ is of similar difficulty
as the problem for general $n$-gons. The only case solved is the one for the octagon settled by
{\sc{Audet, Hansen, Messine}} and {\sc{Perron}} \cite{AudetHansenMessinePerron}. {\sc{Mossinghoff}}
\cite{Mossinghoff08} constructed a sequence of equilateral $n$-gons with unit diameter for $n=2^k$,
$k\ge4$, and proved that their perimeter differs from the maximum perimeter of such $n$-gons by
a term of at most $O(n^{-4})$. By constructing a differen sequence of polygons {\sc{Bingane}} and
{\sc{Audet}} \cite{BinganeAudet21a} further improved the lower bound for the optimum perimeter.

The question about the maximum area of an equilateral $n$-gon with given diameter $d$ is solved for
all $n$: It is $\frac{d^2n}{2}\cos\frac{\pi}{n}\tan\frac{\pi}{2n}$, attained only for a regular $n$-gon.
This follows from Reinhardt's theorem for odd $n$ and was proved by {\sc{Audet}} \cite{Audet17}
for even $n$. {\sc{Bingane}} and {\sc{Audet}} \cite{BinganeAudet21b} determined the equilateral octagon of unit diameter
with maximum width. They also provided a family of equilateral $n$-gons of unit diameter, for $n=2^s$ with
$s\ge4$, whose widths are within $O(n^{-4})$ of the maximum width. It appears that the question about the
maximum width of an equilateral polygon with $n=2^k$ sides and a given perimeter has not been studied so far.

By restricting the class of competing polygons to equilateral polygons some problems with obvious
solutions become interesting. The area, perimeter, and diameter of a general unit-width convex $n$-gon
can be arbitrarily large. This is still the case for an equilateral polygon with an even number of sides.
However, these quantities are bounded for equilateral convex $n$-gons when the number of sides is odd.
{\sc{Audet}} and {\sc{Ninin}} \cite{AudetNinin} determined the maximal perimeter, diameter and area
of an equilateral unit-width convex $n$-gon for every odd $n\ge3$. The optimal polygon is the same for
all three problems: For $n=3$ it is an equilateral triangle of side length $\frac{2}{\sqrt3}$, and
for $n=2k+1\ge5$ a trapezoid whose non-parallel sides have length equal to $\frac{2}{\sqrt3}$,
and the parallel ones have length $m\frac{2}{\sqrt3}$ and $(m-1)\frac{2}{\sqrt3}$.

The papers by {\sc{Mossinghoff}} \cite{Mossinghoff06b} and {\sc{Audet, Hansen}} and {\sc{Messine}}
\cite{AudetHansenMessine07b,AudetHansenMessine09a} contain nice surveys about variations of the
isoperimetric problem for polygons.


\small{
\bibliography{pack}}